\documentclass[12pt]{article}
\usepackage{graphicx,a4}
\usepackage{amsfonts,amsmath,amssymb,mathrsfs,amsthm,epsfig,color}
\usepackage{hyperref}
\definecolor{Red}{cmyk}{0,1,1,0}

\definecolor{verde}{cmyk}{1,0,1,0}

\definecolor{azul}{cmyk}{1,1,0,0}

\usepackage{tocloft}

\newtheorem{theorem}{Theorem}

\newtheorem{corollary}{Corollary}
\newtheorem{lemma}{Lemma}
\newtheorem{proposition}{Proposition}
\newtheorem{theoremalpha}{Theorem}

\theoremstyle{definition}
\newtheorem{definition}[theorem]{Definition}

\newtheorem{example}[theorem]{Example}

\global\newcount\numsec\global\newcount\numfor

\begin{document}
\begin{center}
\vspace*{0.7 cm}
%
%
{\Large\bf
Interactions, Specifications, DLR Probabilities and the
Ruelle Operator in the One-Dimensional Lattice
}
%
%
%
%
\vspace*{0.5cm}\\
    Leandro Cioletti$^{\mbox{\footnotesize{$\dagger$}}}$ and
    Artur O. Lopes$^{\mbox{\footnotesize{$\ddagger$}}}$

\vskip 3mm
$\phantom{.}^{ \mbox{\footnotesize{$\dagger$}}}$
Dep. Matem\'atica - Universidade de Bras\'ilia, 70910-900 Brasilia-DF, Brazil
\\
$\phantom{.}^{ \mbox{\footnotesize{$\ddagger$}}}$
Inst. Matem\'atica UFRGS,  91.500 Porto Alegre-RS, Brazil
\end{center}
\vskip 5mm

%
%
%
%
\begin{abstract}
In this paper, we describe several different meanings for the concept   of Gibbs measure on the
lattice $\mathbb{N}$ in the context of finite alphabets (or state space). We compare and analyze these
``in principle" distinct notions: DLR-Gibbs measures,
Thermodynamic Limit and eigenprobabilities for the dual of the Ruelle operator
(also called conformal measures).

Among other things we extended the classical notion
of a Gibbsian specification on $\mathbb{N}$ in such way that the
similarity of many results in Statistical Mechanics and
Dynamical System becomes apparent. One of our main result claims that the
construction of the conformal Measures in Dynamical
Systems for Walters potentials, using the Ruelle operator,
can be formulated in terms of Specification.
We also describe the Ising model, with $1/r^{2+\varepsilon}$
interaction energy, in the Thermodynamic Formalism setting
and prove that its associated potential is in Walters space - we present an explicit expression.
We also provide an
alternative way for obtaining the uniqueness of the
DLR-Gibbs measures.
\end{abstract}

\section{Introduction}
The basic idea of the Ruelle Operator remounts to the
transfer matrix method introduced by Kramers and Wannier
\cite{MR0004803}
and (independently) by Montroll \cite{montroll1941}, on an effort to compute
the partition function of the Ising model. In a very famous work
published by Lars Onsager in 1944 \cite{MR0010315}, the
transfer matrix method was
generalized to the two-dimensional lattice and was
employed to successfully compute the partition function
for the first neighbors Ising model. As a byproduct, he
obtained the critical point at which the model passes
through a phase transition. These two historical and remarkable
 chapters of the
theory of transfer operators are related to
the study of their actions on finite-dimensional vector spaces.

In a seminal paper in 1968, David Ruelle \cite{MR0234697}
introduced the transfer operator for an one-dimensional
statistical mechanics model with infinite range interactions.
This paved the way to the study of transfer operators in infinite-dimensional
vector spaces. In this paper, Ruelle proved the existence and
uniqueness of the Gibbs measure for a lattice gas system with a potential
depending on infinitely many coordinates.

Nowadays, the transfer operators are
called Ruelle operators (mainly in Thermodynamic Formalism)
and play an important role in Dynamical Systems and
Mathematical Statistical Mechanics.
They are actually useful tools in several other branches of mathematics.

Roughly speaking, the famous Ruelle-Perron-Frobenius Theorem
states that the Ruelle operator
for a potential with a certain regularity,
acting on a suitable Banach space, has a unique simple
positive eigenvalue (equal to the spectral radius)
and associated to it a positive eigenfunction. For H\"older continuous potentials
the proof of this theorem can be found in
\cite{MR1793194,MR1085356,MR0234697}.
In 1978, Walters obtained the Ruelle-Perron-Frobenius Theorem
for a more general setting \cite{MR0466493}, allowing expansive and
 mixing dynamical systems
together with potentials with summable variation.

The Ruelle operator was successfully used to study the problem of existence
and uniqueness of equilibrium states, introduced in \cite{MR0217610,Walters1975},
for a very general class of potentials $f$, see also \cite{MR0404584}.
Under some regularity conditions on $f$ one can show
the uniqueness of the equilibrium states, see
\cite{MR1793194,MR2423393,MR1085356,MR0234697,MR2129258,MR1783787,MR2122435,MR2342978}
and references therein.
Some important properties of the equilibrium probability can be
derived from the  Ruelle operator and this operator
turns out to be a very important tool on topological dynamics
and differentiable dynamical systems,
with applications to the study of invariant measures
for an Anosov diffeomorphism \cite{MR2423393,MR0399421} and the meromorphy of
Zelberg's zeta function \cite{MR1920859}.

The so called  DLR Gibbs measures
were introduced in 1968 and 1969 independently
by Dobrushin \cite{MR0231434}
and Lanford and Ruelle \cite{MR0256687}.
The abstract formulation in terms of specifications
was developed five years latter in\cite{dobrushin1970prescribing,MR0426176,MR0448630}
An important stage in the development of the theory was established by
the works of Preston \cite{MR0448630} and
Gruber, Hintermann and Merlini \cite{MR0523383}  around 1977,
Ruelle (1978) \cite{MR511655} and Israel (1979) \cite{MR517873}.
Preston's work was more focused
on the abstract measure theory, while Gruber et al. concentrated on
specific methods for Ising type models, Israel dealt with
the variational principle and Ruelle worked towards Gibbsian
formalism in Ergodic Theory.

Dobrushin began the study of non-uniqueness of the
DLR Gibbs measure and proposed its interpretation as a
phase transition.
He proved the famous
Dobrushin Uniqueness Theorem in 1968,
ensuring the uniqueness of the Gibbs measures
for a very general class of interactions
at very high temperatures ($\beta \ll 1$).
This result, together with
the rigorous proof of non-uniqueness of
the Gibbs measures for the two-dimensional Ising model
at low temperatures, is a great triumph of the DLR approach
in the study of phase transition in Statistical Mechanics.
Some accounts of the general results on the Gibbs Measure
theory (from the Statistical Mechanics' viewpoint)
can be found in
\cite{MR2252929,MR2189669,MR1241537,Friedli-Velenik17,MR2807681,MR511655}.

In Section 2 of \cite{sarig2009lecture} the author introduces
a concept of DLR-Gibbs measure in the context of topological Markov shifts.
Afterwards the concept and existence of Thermodynamic Limit
were discussed in such context.
Here the definitions of DLR-Gibbs measures and Thermodynamic Limit
are similar to the ones considered in \cite{sarig2009lecture}.
We shall remark that in reference \cite{sarig2009lecture} (see Definition 1.4)
the concept of Gibbs measure is considered in the sense of Bowen.
Here we will not work with this concept of Gibbs measure.

The present work aims to explain how to use
DLR-Gibbs measures to obtain the conformal measures
considered in Thermodynamic Formalism.
In order to do that we introduce a notion of specification associated
to continuous potential. In particular, we
show how to construct an absolutely uniformly summable
specification for any H\"older potential and use this construction
to motivate the specifications considered here.
The main results of this paper are
Theorems \ref{Teo-DLR-Lim-Term} and \ref{teo-principal}, in Section \ref{sec-main-results},
which prove the equivalence between the conformal measures
considered in Thermodynamic Formalism and DLR-Gibbs measures,
for potentials in the Walters space.

\medskip

The Preprint \cite{CL-rcontinuas-2016} approaches similar problems (as described here)
but in different setting. For example, potentials can be continuous functions
and the alphabet can be any compact metric space (which includes uncoutable alphabets).
But, on the other hand, the strong equivalence proved here in Theorem \ref{teo-principal}
is no longer true in this setting.

\section{Ruelle Operator and Conformal Measures}
In this paper $\mathbb{N}$ denotes the set of positive integers,
$\mathscr{A}$ is a finite alphabet and
$\Omega\equiv \mathscr{A}^\mathbb{N}$ denotes the symbolic
space endowed with its standard
metric $d$ given by
$d(x,y) = 2^{-N}$, where $N = \inf \{i\in\mathbb{N} : x_i\neq y_i\}$.
The Borel $\sigma$-algebra of $\Omega$ is denoted by $\mathscr{F}$.
The dynamics here is given by $\sigma:\Omega\to\Omega$, the left-shift mapping.
The space of all real continuous bounded functions on $\Omega$ endowed with its
standard supremum norm $\|\cdot\|_{\infty}$
is denoted simply by $C(\Omega)$. We use the notation
$\mathcal{P}(\Omega)\equiv \{ \nu:\mathscr{F}\to [0,1]:\ \nu\ \text{is a probability measure} \}$
for the set of all Borel probability measures over $\Omega$.
\begin{definition}[Ruelle Operator]
Let $f:\Omega\to\mathbb{R}$ be a continuous function.
The Ruelle operator associated to $f$, notation
$\mathcal{L}_f:C(\Omega)\to C(\Omega)$
is defined on the function $\psi$ as follows
\[
\mathcal{L}_{f} (\psi)(x)
=
\sum_{y\in\Omega;\ \sigma(y)=x} \exp(f(y)) \, \psi(y).
\]
\end{definition}
Normally we call $f$ a potential and $\mathcal{L}_{f}$
the transfer operator associated to the potential $f$.
The dual of the Ruelle operator $\mathcal{L}_{f}^{*}$  acts on the set of
Borel finite signed measures over $\Omega$ as follows
$\mathcal{L}^{*}_{f}(\nu)(\psi)=  \nu(\mathcal{L}_{f}(\psi))$
for all $\psi\in C(\Omega)$.
	
Fix $0<\alpha<1$. We say that a function $f:\Omega\to\mathbb{R}$ is
$\alpha$-H\"older continuous if
\[
\mathrm{Hol}_{\alpha}(f)\equiv \sup_{x\neq y}\frac{|f(x)-f(y)|}{d^{\alpha}(x,y)} <+\infty.
\]
The space of all real $\alpha$-H\"older continuous functions on $\Omega$ is denote
by $C^{\alpha}(\Omega)$. When we say that $f$
is H\"older continuous function we mean $f\in C^{\alpha}(\Omega)$
for some $0<\alpha<1$.
 For any $n\geq 1$ we define
the $n$-th variation of a function $f:\Omega \to \mathbb{R}$ by
$
\mathrm{var}_n (f) =
\sup \, \{\,| f(x)-f(y) |: x, y\in  \Omega \,\,\text{and}\
x_i =y_i\, \, \text{ for all}\,\,  1\leq  i\leq n\}.
$
We say that a function $f:\Omega \to \mathbb{R}$ is in the {\bf Walters} space, notation $W(\Omega)$,
if the following condition is satisfied
\begin{align}\label{walters-condition}
\lim_{p\to \infty}\,
\,\sup_{n\geq 1}
\,\mathrm{var}_{n+p} (S_n(f))	= 0,
\quad \text{where}\ \ S_n(f)=f+\ldots+f\circ \sigma^{n-1}.
\end{align}
We remark that for any $0<\alpha<1$ we have
$
C^{\alpha}(\Omega)\subset W(\Omega)\subset C(\Omega).
$

\begin{theorem}[Ruelle-Perron-Frobenius (RPF) for Walters Potentials]
\label{RPF}
	Let $f$ be a potential in $W(\Omega)$.
	Then there exists a strictly positive function $\psi_{f}\in W(\Omega)$
	and a strictly positive eigenvalue $\lambda_{f}$ such
	that $\mathcal{L}_{f}(\psi_f)=\lambda_{f}\psi_{f}$.
	The eigenvalue $\lambda_{f}$ is simple and it is equal
	to the spectral radius of the operator.
	Moreover, there exists a unique probability measure
	$\nu_f$ over $\Omega$ such that
	$\mathcal{L}_{f}^*(\nu_f)=\lambda_f\, \nu_f$.
\end{theorem}

\begin{proof}
For a proof see \cite{MR1841880,MR1783787,MR2122435,MR2342978}.
\end{proof}

\begin{definition}
Let $f\in C(\Omega)$ a continuous potential and $\rho(\mathcal{L}_{f})$
the spectral radius of $\mathcal{L}_{f}$ acting on $C(\Omega)$.
 The set of all Borel probability measures $\nu$ over $\Omega$,
 satisfying $\mathcal{L}^{*}_{f}\nu = \rho(\mathcal{L}_{f})\nu $
 is denoted by $\mathcal{G}^{*}(f)$.
\end{definition}

Note that if $f\in W(\Omega)$,
then follows from Theorem \ref{RPF} that $\rho(\mathcal{L}_{f}) = \lambda_{f}$
and
\[
\mathcal{G}^{*}(f)
=\{\nu\in \mathcal{P}(\Omega): \mathcal{L}^{*}_{f}\nu = \lambda_{f}\nu  \}
\]
is a singleton.

\section{Interactions and Continuous Potentials} \label{Spec}

In the classical literature on Statistical Mechanics the concept of interaction is prominent.
In what follows we described it but only in the generality needed in this paper.
For a comprehensive exposition on this topic see \cite{MR2807681}.

From now on the notation $A\Subset \mathbb{N}$
means that
$A$ is an empty or finite subset of $\mathbb{N}$.
If for each $A\Subset \mathbb{N}$ we
associated a function $\Phi_{A}:\Omega\to\mathbb{R}$
then we have a family of functions defined on $\Omega$
and indexed on the finite parts of $\mathbb{N}$.
We denote such family simply by $\Phi=\{\Phi_A\}_{A\Subset \mathbb{N}}$
and $\Phi$ will be called an interaction.
We shall remark that is usual $\Phi$
to have several finite subsets $A$'s for which the
associated function $\Phi_{A}$ is identically zero.

The space of interactions has natural structure of a vector space
where the sum of two interactions $\Phi$ and $\Psi$,
is given by the interaction $(\Phi+\Psi)\equiv\{\Phi_{A}+\Psi_{A}\}_{A\Subset \mathbb{N}}$
and $\lambda\Phi=\{\lambda \Phi_{A}\}_{A\Subset \mathbb{N}}$,
for any $\lambda\in\mathbb{R}$.
This vector space is too big for our purposes so
we focus in a proper subspace of it.

Before proceed we shall remark that
we can also consider interactions defined
on a general countable set $V$.
If $V=\mathbb{Z}$, for example,
then the family $\Phi$ is now indexed over the
collection of all $A\Subset \mathbb{Z}$.
In this case we say that the interaction is
defined on the lattice $\mathbb{Z}$.
We focus here on interactions $\Phi$
defined on the lattice $\mathbb{N}$, in order to
relate the DLR-Gibbs measures and
the Thermodynamic Formalism.

\begin{definition}[Uniformly Absolutely Summable Interaction]
An interaction $\Phi=\{\Phi_{A}\}_{A\Subset\mathbb{N}}$ is
called {\it uniformly absolutely summable} (UAS) interaction if it satisfies:
\begin{enumerate}
\item
for each $A\Subset \mathbb{N}$ the function $\Phi_{A}:\Omega\to\mathbb{R}$
depends only on the coordinates with indexes in $A$;
\item
$
\displaystyle
		\|\Phi\|\equiv		
		\sup_{n\in\mathbb{N}}
		\sum_{ \substack{A\Subset \mathbb{N}; A\ni n} }
		\sup_{x\in\Omega}|\Phi_{A}(x)|<\infty.
$
\end{enumerate}
\end{definition}
\begin{example}[Dyson Model on $\mathbb{N}$]
\label{exemplo-ising-longo-alcance}
Consider the alphabet $\mathscr{A}=\{-1,1\}$
and a fixed $\alpha>1$. Then the interaction $\Phi$ given by
	\[
		\Phi_{A}(x) =
		\begin{cases}
			\displaystyle\frac{x_{n}x_{m}}{|n-m|^{\alpha}},
					&\ \text{if}\ \ A=\{n,m\}\ \text{and}\ m\neq n;
			\\
			0,			&\ \text{otherwise},
		\end{cases}
	\]
is an UAS interaction.
In fact, for any $A\Subset \mathbb{N}$ we have that
$\Phi_{A}\equiv 0$ if $\#A\neq 2$.
On the other hand, if $A=\{m,n\}$ with $m\neq n$
we have that $\Phi_{A}$ depends only on the coordinates $x_{n}$ and $x_{m}$.
The regularity condition is verified as follows
\begin{align*}
		\|\Phi\|\equiv		
		\sup_{n\in\mathbb{N}}
		\sum_{ \substack{A\Subset \mathbb{N}; A\ni n} }
		\sup_{x\in\Omega}|\Phi_{A}(x)|
		&=
		\sup_{n\in\mathbb{N}}
		\,
		\sum_{ \substack{m\in \mathbb{N}\setminus\{n\}} }
		\sup_{x\in\Omega} \frac{|x_nx_m|}{|m-n|^{\alpha}}
		\leq
		2\zeta(\alpha),
\end{align*}
where $\zeta$ is the Riemann zeta function.
\end{example}
In order to state our next result we introduce some notations.
Let $x,y$ and $z\in \Omega$ and $n,m\mathbb{N}$. We use the notation
$x_{1}^{n}y_{n+1}^{n+m}z_{n+m+1}^{\infty}$ to denote a point in $\Omega$,
whose its coordinates are $(x_1,\ldots,x_n,y_{n+1},\ldots,y_{n+m},z_{n+m+1},\ldots)$.
For each $k\geq 1$ and $n\geq 0$ consider the arithmetic progression
$A(k,n)\equiv\{k,\ldots,2k+n\}$. For each $f\in C(\Omega)$ and $y\in\Omega$
we define $f_{A(k,n)}:\Omega\to\mathbb{R}$ as follows
\[
		f_{A(k,n)}(x)
		=
		f(x_k^{2k+n}y_{2k+n+1}^{\infty})
		-
		f(x_k^{2k+n-1}y_{2k+n}^{\infty})
\]
if $n\geq 1$ and
$		
f_{A(k,0)}(x)
=
f(x_{k}^{2k}y_{2k+1}^{\infty})
-		
f(y).
$

\begin{lemma}\label{lema-fcont-interacao}
Let $f\in C(\Omega)$, $y\in\Omega$ and
$\Phi^{f}\equiv \{\Phi^{f}_{A}\}_{A\Subset \mathbb{N}}$ be the interaction
given by
$\Phi^{f}_{A}(x)=f_{A(k,n)}(x)$ if $A=A(k,n)$ and $0$ otherwise.
Then for all $x\in\Omega$ we have
\[
		f(x)=
		f(y)+\sum_{ \substack{A\Subset \mathbb{N}; A\ni 1} }
		\Phi_{A}^{f}(x).
\]
\end{lemma}
\begin{proof}
For any $n\geq 0$ we have
$
\sum_{j=0}^{n}
\Phi^{f}_{A(1,j)}(x)
=
(f(x_1^2y_3^{\infty})-f(y))
+(f(x_1^{3}y_4^{\infty})-f(x_1^{2}y_3^{\infty}))+
\ldots +
	(f(x_1^{n+2}y_{n+3}^{\infty})
	-f(x_1^{n+1}y_{n+2}^{\infty}))
=
f(x_1^{n+2}y_{n+3}^{\infty})-f(y).
$
From the continuity of $f$ and the previous equation the lemma follows.
\end{proof}

\begin{proposition}\label{prop-holder-AUS}
Let $0<\alpha<1$, $f\in C^{\alpha}(\Omega)$ and $\Phi^{f}$ as
in previous lemma. Then $\Phi^{f}$ is an UAS interaction.
\end{proposition}
\begin{proof}
Note that for all $n\geq 1$ we have
\[
\sum_{ \substack{A\Subset \mathbb{N}; A\ni n} }
		\|\Phi^{f}_{A}(x)\|_{\infty}
\leq
\sum_{k=1}^n\sum_{m=k+1}^{\infty}\|\Phi^{f}_{A(k,m)}(x)\|_{\infty}
\]	
Since $f\in C^{\alpha}(\Omega)$, for all $k\in \{1,\ldots,n\}$ and $m\geq k+1$ we have
$
\|\Phi^{f}_{A(k,m)}(x)\|_{\infty}
\leq
\mathrm{Hol}_{\alpha}(f) 2^{-\alpha (k+m-1)}.
$
Therefore
	\begin{align*}
		\|\Phi^{f}\|
		&\equiv		
		\sup_{n\in\mathbb{N}}
		\sum_{ \substack{A\Subset \mathbb{N}; A\ni n} }
		\|\Phi^{f}_{A}\|_{\infty}
		\leq
		\frac{2^{\alpha}\mathrm{Hol}_{\alpha}(f)}{(1-2^{-\alpha})^2}.
\qedhere
\end{align*}	
\end{proof}

Let $\Phi$ be a UAS interaction. For each $n\in\mathbb{N}$ let $\Lambda_n\equiv \{1,\ldots,n\}$.
The function
	\begin{equation}\label{hamiltoniano-geral}
		H_{n}(x)
		=
		\sum_{\substack{ A\Subset \mathbb{N}
						\\
						A\cap \Lambda_n\neq \emptyset} }
		\Phi_A(x)
	\end{equation}
is called the Hamiltonian associated to the interaction $\Phi$
in the volume $\Lambda_n$.
In Mathematical Statistical Mechanics the Gibbs measures
(called here DLR-Gibbs measures) associated to an interaction
is normally constructed by means of $(H_n)_{n\geq 1}$.
Before explain this construction we obtain a formula for $H_n$
when $\Phi\equiv\Phi^{f}$ is a UAS interaction.

\begin{proposition}\label{propo-relacao-hamiltoniano-potencial}
Let $f\in C(\Omega)$ and assume that $\Phi^{f}$ defined as in
Lemma \ref{lema-fcont-interacao} is a UAS interaction.
Then, there is a constant $C$ so that for all $n\in\mathbb{N}$
the Hamiltonian $H_n$ defined by \eqref{hamiltoniano-geral}
satisfies
\[
H_n(x)
=
f(x)+\ldots+f(\sigma^{n-1} x)
+nC.
\]
\end{proposition}
\begin{proof}
From the definition of $\Phi^f$ and the UAS property we have
	\begin{align*}
		H_{n}(x)
		=
		\sum_{\substack{ A\Subset \mathbb{N}
						\\
						A\cap \Lambda_n\neq \emptyset} }
		\Phi^{f}_A(x)
		=
		\sum_{k=1}^{n}
		\sum_{m=0}^{\infty}
		\Phi_{A(k,m)}(x).
	\end{align*}
By using similar argument as in Lemma \ref{lema-fcont-interacao}
we can prove that the inner sum in rhs above is given by
$
f(\sigma^{k-1} x)-f(y).
$
By taking $C\equiv f(y)$ and them summing the last
expression with $k$ varying from $1$ to $n$ the proposition follows.
\end{proof}

Aiming to have an equivalent description of
conformal measures associated to a Walters potential $f$
and the DLR-Gibbs measure associated to $\Phi^f$
we develop below the theory of DLR-Gibbs measures (within our setting)
dispensing the UAS hypothesis.

\section{Specifications and DLR-Gibbs Measures}

From now on, the Hamiltonian $H_n$ is assumed to be of form
\begin{align}\label{Hn-f-cont}
H_n(x)
=
f(x)+f(\sigma x)+\ldots+f(\sigma^{n-1} x)
+nC,
\end{align}
where $f\in C(\Omega)$. In this section we extend some classical
results about Gibbsian specifications to the case where the Hamiltonian has
the above form. The motivation to extend the DLR theory on this direction
becomes natural in view of the results of the previous section and
this extension is crucial to show the equivalence stated in Theorem \ref{teo-principal}.

\begin{lemma}\label{lema-Hn-interchange}
For any $n,r\in\mathbb{N}$, $x,y$ and $z\in\Omega$ we have
	\[
	H_{n+r}(y_1^{n+r}z_{n+r+1}^{\infty})
	\!-\!
	H_{n}(y_1^{n+r}z_{n+r+1}^{\infty})
	\!=\!
	H_{n+r}(x_{1}^{n}y_{n+1}^{n+r}z_{n+r+1}^{\infty})
	\!-\!
	H_{n}(x_1^{n}y_{n+1}^{n+r}z_{n+r+1}^{\infty}).
	\]
\end{lemma}
\begin{proof}
From definition of $H_n$ follows that
\[
	H_{n+r}(y_1^{n+r}z_{n+r+1}^{\infty})
	-
	H_{n}(y_1^{n+r}z_{n+r+1}^{\infty})
	=
	\sum_{j=n}^{n+r-1}f(\sigma^j(y_1^{n+r}z_{n+r+1}^{\infty})).
\]
Since rhs above is equals to
$
H_{n+r}(x_{1}^{n}y_{n+1}^{n+r}z_{n+r+1}^{\infty})
\!-\!
H_{n}(x_1^{n}y_{n+1}^{n+r}z_{n+r+1}^{\infty})
$ the lemma is proved.
\end{proof}

\begin{definition}\label{def-specificacao-determinada-Phi}
Given a continuous potential $f$  we define a family of
probability kernels $(K_n)_{n\geq 1}$, where
for each $n\in\mathbb{N}$ the kernel $K_n:\mathscr{F}\times\Omega\to\mathbb{R}$
is given by
\[
K_{n}(F,y) = \frac{1}{Z_{n}^{y}}
		\sum_{ \substack{ x\in\Omega;\\ \sigma^n(x)=\sigma^n(y) }   }
\!\!\!\!1_{F}(x)\exp(H_n(x)),
\ \text{where}\ \
Z_{n}^{y}
=
\!\!\!\!\!
\sum_{ \substack{ x\in\Omega;\\ \sigma^n(x)=\sigma^n(y) }   }
\exp(H_n(x)).
\]
\end{definition}

Note that the constant $C$ in \eqref{Hn-f-cont} is irrelevant for
the definition of $K_n$, therefore, without loss of generality,
we can assume that $C= 0$.
Let $\pi_n:\Omega\to \mathscr{A}$ be the canonical projection in the $n$-th
coordinate and $\mathcal{T}_n$ the sigma-algebra generated by
the projections $\{\pi_j: j\geq n+1\}$.
Then for any $f\in C(\Omega)$ and for all $n\in \mathbb{N}$
it is easy to see that the kernel $K_n$ satisfies:
	\begin{enumerate}
		\item[a)] $y\longmapsto K_{n}(F,y)$ is $\mathcal{T}_n$-measurable;
		\item[b)] $F\longmapsto K_{n}(F,y)$ is a Borel probability measure;
		\item[c)] $y\longmapsto \int_{\Omega}g(x)\, dK_{n}(x,y)$ is continuous for any $g\in C(\Omega)$.
	\end{enumerate}

\begin{theorem}[Compatibility Conditions]\label{teo-DLR-finite-vol}
	If $(K_n)_{n\geq 1}$ is a family of probability kernels as in
	Definition \ref{def-specificacao-determinada-Phi}, then
    for each fixed $z\in\Omega$ and for any integers $r,n\geq 1$
    we have
	\[
		\int_{\Omega}
			\left[
				\int_{\Omega} g(x) \, dK_{n}(x,y)
			\right]
			\, dK_{n+r}(y,z)
		=
		\int_{\Omega} g(y) \, dK_{n+r}(y,z),
		\qquad \forall\ g\in C(\Omega).
	\]
\end{theorem}
\begin{proof}
Follows from the definition of $K_n$ that for any $g\in C(\Omega)$
\begin{align*}
	\int_{\Omega} g(x) \, dK_{n}(x,(y_1^{n+r}z_{n+r+1}^{\infty}))
	&=
	\frac{1}{Z_{n}^{(y_{1}^{n+r}z_{n+r+1}^{\infty})}}
	\sum_{\substack{ x\in\Omega  \\ \sigma^n(x)=(y_{n+1}^{n+r}z_{n+r+1}^{\infty})   }  }
	\!\!\!\!\!\!\!
	g(x) \exp(H_n(x))
	\\
	&\equiv
	h(y_{1}^{n+r}z_{n+r+1}^{\infty}).
\end{align*}
We are using above the notation $h(y_{1}^{n+r}z_{n+r+1}^{\infty})$ for the sake of compatibility, but
note that this quantity do not depends on $y_1,\ldots,y_n$.

Therefore to prove the theorem is enough to show that
	\[
		\frac{1}{Z_{n+r}^{z}}	
		\!\!\!\!\!\!\!\!
		\sum_{ \substack{ y\in\Omega \\ \sigma^{n+r}(y)=\sigma^{n+r}(z) }  }
		\!\!\!\!\!\!\!\!\!\!
		h(y)
		\exp(H_{n+r}(y))
		=
		\frac{1}{Z_{n+r}^{z}}	
		\!\!\!\!\!\!\!\!
		\sum_{ \substack{ y\in\Omega \\ \sigma^{n+r}(y)=\sigma^{n+r}(z) }  }
		\!\!\!\!\!\!\!\!\!\!
		g(y)
		\exp(H_{n+r}(y)).
	\]
Since $Z_{n+r}^{z}>0$, the above equation is equivalent to
	\begin{equation}\label{eq-compatibilidade-gibbs-1}
		\sum_{ \substack{ y\in\Omega \\ \sigma^{n+r}(y)=\sigma^{n+r}(z) }  }
		\!\!\!\!\!\!\!\!\!\!
		h(y)
		\exp(H_{n+r}(y))
		=
		\sum_{ \substack{ y\in\Omega \\ \sigma^{n+r}(y)=\sigma^{n+r}(z) }  }
		\!\!\!\!\!\!\!\!\!\!
		g(y)
		\exp(H_{n+r}(y)).
	\end{equation}
In order to prove the theorem we show in the sequel
that \eqref{eq-compatibilidade-gibbs-1} holds.
Indeed, from the definition of $h$, we have that the l.h.s above is given by
	\[
		\sum_{ \substack{ y\in\Omega \\ \sigma^{n+r}(y)=\sigma^{n+r}(z) }  }
		\!
		\frac{1}{Z_{n}^{(y_1^{n+r}z_{n+r+1}^{\infty})}}
		\sum_{\substack{ x\in\Omega  \\ \sigma^n(x)=(y_{n+1}^{n+r}z_{n+r+1}^{\infty})}  }
		g(x) \exp(H_n(x)+H_{n+r}(y)).
	\]
From Lemma \ref{lema-Hn-interchange} follows that the above expression
is equal to
	\[
		\sum_{ \substack{ y\in\Omega \\ \sigma^{n+r}(y)=\sigma^{n+r}(z) }  }
		\!\!\!\!\!\!
		\frac{\exp(H_n(y_1^{n+r}z_{n+r+1}^{\infty}))}{Z_{n}^{(y_1^{n+r}z_{n+r+1}^{\infty})}}
		\sum_{\substack{ x\in\Omega  \\ \sigma^n(x)=(y_{n+1}^{n+r}z_{n+r+1}^{\infty}) }  }
		\!\!\!\!\!\!\!\!\!\!\!\!
		g(x)
		\exp(H_{n+r}(x_1^{n}y_{n+1}^{n+r}z_{n+r+1}^{\infty})).
	\]
Note that the above expression is equals to
	\[
		\sum_{ y_1,\ldots,y_{n+r}\in\mathscr{A} }
		\frac{\exp(H_n(y_1^{n+r}z_{n+r+1}^{\infty}))}{Z_{n}^{(y_{1}^{n+r}z_{n+r+1}^{\infty})}}
		\sum_{x_1,\ldots,x_n\in\mathscr{A} }
		g(x)
		\exp(H_{n+r}(x_1^{n}y_{n+1}^{n+r}z_{n+r+1}^{\infty})).
	\]
By interchanging summation order we can rewrite the above expression as
	\[
		\sum_{ y_{n+1},\ldots,y_{n+r}\in\mathscr{A} }
		\sum_{ y_1,\ldots,y_{n}\in\mathscr{A} }
		\!\!\!\!
		\frac{\exp(H_n(y_1^{n+r}z_{n+r+1}^{\infty}))}{Z_{n}^{(y_{1}^{n+r}z_{n+r+1}^{\infty})}}
		\!\!\!\sum_{x_1,\ldots,x_n\in\mathscr{A} }
		\!\!\!\!
		g(x)
		\exp(H_{n+r}(x_1^{n}y_{n+1}^{n+r}z_{n+r+1}^{\infty})).
	\]
Since the third sum above do not depend on $y_1,\ldots,y_n$ and
\[
\sum_{ y_1,\ldots,y_{n}\in\mathscr{A} }
\!\!\!\!
\frac{\exp(H_n(y_1^{n+r}z_{n+r+1}^{\infty}))}{Z_{n}^{(y_{1}^{n+r}z_{n+r+1}^{\infty})}}
=1
\]
the previous expression is equal to
	\[
		\sum_{ \substack{  y_{n+1},\ldots,y_{n+r}\in\mathscr{A} \\ x_1,\ldots,x_n \in\mathscr{A}   }  }
		g(x)
		\exp(H_{n+r}(x_1^{n}y_{n+1}^{n+r}z_{n+r+1}^{\infty}))
		=
		\sum_{ \substack{ y\in\Omega \\ \sigma^{n+r}(y)=\sigma^{n+r}(z) }  }
		\!\!\!\!\!\!\!\!\!\!
		g(y)
		\exp(H_{n+r}(y)).
	\]
The last expression shows
that \eqref{eq-compatibilidade-gibbs-1} holds
and the theorem is proved.
\end{proof}

Notice that the collection $(K_n)_{n\in\mathbb{N}}$ is similar to but not exactly
a quasilocal specification as in the literature of Mathematical Statistical Mechanics,
see for example \cite{MR1241537,MR2807681,MR0448630}. It is possible to extend this
collection to a classical quasilocal specification, but the point here is to obtain
similar results to the classical theory of DLR-Gibbs measures in this more general
setting. For the extension argument see \cite{CL-rcontinuas-2016}.

\begin{proposition}\label{prop-DLR-equations}
	Let $(K_n)_{n\geq 1}$ be as in Definition \ref{def-specificacao-determinada-Phi} and	
	$z\in \Omega$ a fixed point.
	If the sequence $K_{n_j}(\cdot,z)\rightharpoonup \mu^{z}$ (weak-$*$ topology),
	when $j\to\infty$,
	then for any continuous function $g:\Omega\to\mathbb{R}$, we have
\[
	\int_{\Omega}
		\left[
			\int_{\Omega} g(x) \, dK_{n}(x,y)
		\right]
	\, d\mu^z(y)
	=
	\int_{\Omega}g\, d\mu^z.
\]
\end{proposition}
\begin{proof}
For any fixed $n\in\mathbb{N}$, the mapping
\[
\Omega\ni y\longmapsto \int_{\Omega} g(x) \, dK_{n}(x,y)
\]
is continuous. From the compatibility condition and the
definition of weak-$*$ topology follows that
\begin{align*}
	\int_{\Omega}
		\left[
			\int_{\Omega} g(x) \, dK_{n}(x,y)
		\right]
	\, d\mu^z(y)
	&=
	\lim_{j\to\infty}
	\int_{\Omega}
		\left[
			\int_{\Omega} g(x) \, dK_{n}(x,y)
		\right]
	\, dK_{n_j}(y,z)
	\\
	&=
		\lim_{j\to\infty}
		\int_{\Omega}
			 g(y) 			
		\, dK_{n_j}(y,z)
	\\
	&=
	\int_{\Omega}g\, d\mu^z.
\qedhere
\end{align*}
\end{proof}

\begin{definition}[DLR-Gibbs Measures]
\label{definicao-gibbs-measures}
Let $(K_n)_{n\in\mathbb{N}}$ be as in Definition \ref{def-specificacao-determinada-Phi}.
The set of DLR Gibbs measures associated to a continuous potential $f$ is defined as
\[
\mathcal{G}^{DLR}(f)
\equiv
	\left\{
		\mu\in \mathcal{P}(\Omega):
		\begin{array}{c}
		\mu(F|\mathcal{T}_n)(y)=K_{n}(F,y)
		\ \text{for}\ \mu-\text{a.a.}\ y,\
		\\
		\forall F\in\mathscr{F}
		\ \text{and}\ \forall n\in\mathbb{N}
		\end{array}
	\right\}.
\]
\end{definition}

The DLR equations play an important role in Statistical Mechanics.

\begin{theorem}[DLR-equations]\label{teo-DLR-equation-equiv-gibbs}
Let $(K_n)_{n\in\mathbb{N}}$ be as in Definition \ref{def-specificacao-determinada-Phi}.
A Borel probability measure $\mu\in\mathcal{P}(\Omega)$ belongs to $\mathcal{G}^{DLR}(f)$
iff for all $n\in\mathbb{N}$ and any
continuous function $g:\Omega\to\mathbb{R}$,
we have
\[
	\int_{\Omega}
		\left[
			\int_{\Omega} g(x) \, dK_{n}(x,y)
		\right]
	\, d\mu(y)
	=
	\int_{\Omega}g\, d\mu.
\]
\end{theorem}
\begin{proof}
We follow closely the reference \cite{MR2807681}.
Suppose that $\mu\in\mathcal{G}^{DLR}(f)$ then
it follows from the definition of $\mathcal{G}^{DLR}(f)$
and the basic properties of the conditional
expectation that for all $n\in\mathbb{N}$ we have
\[
	\int_{\Omega} g\, d\mu
	=	
	\int_{\Omega} \mu(g|\mathcal{T}_n)(y) \, d\mu(y)
	=
	\int_{\Omega}
		\left[
			\int_{\Omega} g(x)\, dK_{n}(x,y)
		\right]
		\, d\mu(y).
\]
Conversely, we assume that the DLR-equations are valid for all $n\in\mathbb{N}$
and for any continuous function $g$.
Let $g=1_{E}h$, where $E\in\mathcal{T}_n$
is a cylinder set and $h$ is an arbitrary continuous function.
Then $g$ is continuous and
\[
	\int_{\Omega}\!\! 1_{E}(y)\!
		\left[
			\int_{\Omega} h(x)\, dK_{n}(x,y)
		\right]
		\!\!d\mu(y)
	\!=\!
	\int_{\Omega}
		\left[
			\int_{\Omega} 1_{E}(x)h(x)\, dK_{n}(x,y)
		\right]
		\!\!d\mu(y)
	=
	\int_{E} h\, d\mu,
\]
where in the first equality we used that the
function $1_{E}$ do not depends on its $n$ first coordinates and
definition of $K_n(\cdot,y)$.
From the  Dominate Convergence Theorem follows that the class of $E$'s satisfying the
above identity is a monotone class, and from
Monotone Class Theorem follows that the above
identity holds for any measurable set $E\in \mathcal{T}_n$.
Since the mapping
\[
	y\mapsto \int_{\Omega} h(x)\, dK_{n}(x,y)
\]
is  $\mathcal{T}_n$-measurable and $E\in\mathcal{T}_n$ is
an arbitrary measurable set,
we have, from the definition of conditional expectation and
last equality, that
\[
	\int_{\Omega} h(x)\, dK_{n}(x,y)
	=
	\mu(h|\mathcal{T}_n)(y)
	\ \ \
	\mu\ \text{a.e.}
\]
Using again
the Dominate Convergence Theorem for conditional expectation and
Monotone Class Theorem
we can show that the above equality holds for $h=1_{F}$
where $F$ is a measurable set in $\mathscr{F}$,
so the result follows.
\end{proof}

From item c) that appears before Theorem \ref{teo-DLR-finite-vol} and
DLR-equations follows that $\mathcal{G}^{DLR}(f)$ is a closed subset
of $\mathcal{P}(\Omega)$, with respect to the weak-$*$ topology.
Since $\mathcal{P}(\Omega)$ endowed with this topology is a compact Hausdorff
space follows that $\mathcal{G}^{DLR}(f)$ is compact.

\medskip

Let $f\in C(\Omega)$ and $(K_n)_{n\in\mathbb{N}}$ as
in Definition \ref{def-specificacao-determinada-Phi}.
For each $y\in\Omega$ we define $\mathscr{C}_y$ as
being the set of all the cluster points, in the weak-$*$ topology,
of the set $\{K_n(\cdot,y): n\geq 1\}$. We call $\mu\in\mathscr{C}_y$
a Thermodynamic Limit obtained from the boundary condition $y$.

\begin{definition}
The closure, in the weak-$*$ topology, of the convex hull of
the set $\cup_{y\in\Omega} \mathscr{C}_y$ will be denoted
by $\mathcal{G}^{TL}(f)$.
\end{definition}

\begin{proposition}\label{prop-GTL-in-GDLR}
For any $f\in C(\Omega)$ we have that
the set $\mathcal{G}^{TL}(f)$ is always non-empty. Moreover,
$\mathcal{G}^{TL}(f)\subset \mathcal{G}^{DLR}(f)$.
\end{proposition}

\begin{proof}
For any compact metric space $\Omega$ we have that $\mathcal{P}(\Omega)$
is compact, with respect to the weak-$*$. Since this topology is metrizable
follows that  $\mathcal{P}(\Omega)$ is sequentially compact.
Therefore the subset $\{K_n(\cdot,y): n\geq 1\}\subset \mathcal{P}(\Omega)$
has at least one cluster point
$\mu^y$, thus proving that $\mathcal{G}^{TL}(f)\neq \emptyset$.
The inclusion
is straightforward application of Proposition \ref{prop-DLR-equations}
and Theorem \ref{teo-DLR-equation-equiv-gibbs}.
\end{proof}

Examples where  one can get different Thermodynamic Limits $\mu\in\mathscr{C}_y$
depending of  the boundary condition $y$ appear in \cite{MR3350377}.

\section{Specifications and Ruelle Operator}

In this section we establish relevant relations, in this work, between
the kernels $(K_n)_{n\in\mathbb{N}}$
given by Definition \ref{def-specificacao-determinada-Phi}
and Ruelle operator $\mathcal{L}_{f}$.

We first recall that the $n$-th iterated of Ruelle operator applied
to any $\psi\in C(\Omega)$ and calculated at $y$
is given by the following formula
\[
\mathcal{L}_{f}^n (\psi)(y)
=
\sum_{ \substack{ x\in\Omega;\\ \sigma^n(x)=y	 }}
\exp(S_n(f)(x))\psi(x).
\]

\begin{proposition}\label{prop-Kn-Ln}
Let $f\in C(\Omega)$. For any cylinder set $F\in\mathscr{F}$
and $n\in\mathbb{N}$ we have
\begin{align*}
		K_n(F,y) = \frac{1}{Z_{n}^{y}}
		\sum_{ \substack{ x\in\Omega;\\ \sigma^n(x)=\sigma^n (y) }   }
		1_{F}(x)\exp(H_n(x))
		=
		\frac{\mathcal{L}_{f}^n (1_F )(\sigma^n (y))}{ \mathcal{L}_{f}^n(1 )(\sigma^n (y)) }.
\end{align*}
\end{proposition}

\begin{proof}
The first equality is simply definition of $K_n$.
From definition we have $H_n(x)=S_n(f)(x)$ so the second equality above
follows from the formula for the $n$-th iterated of Ruelle operator
since
\begin{align*}
\mathcal{L}_{f}^n (1_F)(\sigma^n(y))
=
\sum_{ \substack{ x\in\Omega;\\ \sigma^n(x)=\sigma^n(y) }}
\exp(S_n(f)(x))1_{F}(x)
=
\sum_{ \substack{ x\in\Omega;\\ \sigma^n(x)=\sigma^n (y) }   }
1_{F}(x)\exp(H_n(x))
\end{align*}
and
\[
\mathcal{L}_{f}^n (1)(\sigma^n(y))
=
\sum_{ \substack{ x\in\Omega;\\ \sigma^n(x)=\sigma^n(y) }}
\exp(S_n(f)(x))
=
Z_n^y.
\qedhere
\]
\end{proof}

\begin{lemma}\label{lemma-comp-op-ruelle}
Let $f$ be a continuous potential.
For all $n,m\in\mathbb{N}$,  $z\in\Omega$ and $\psi\in C(\Omega)$ we have
\[
  \mathcal{L}^{n+m}_f (\psi) (\sigma^{n+m} (z))
  =
  \mathcal{L}^{n+m}_f\,
  	\left(
  		\frac{\mathcal{L}^{n}_f(\psi)(\sigma^n (\cdot))}
  			{\mathcal{L}^{n}_f(1)(\sigma^n (\cdot))}
  	\right)\!\!
    (\sigma^{n+m} (z)).
\]
\end{lemma}
\begin{proof}
The proof follows from Proposition \ref{prop-Kn-Ln} and
Theorem \ref{teo-DLR-finite-vol} (compatibility conditions for $(K_n)_{n\in\mathbb{N}}$).
\end{proof}

\section{Main Results}\label{sec-main-results}

\begin{lemma}\label{lemma-aux-unicidade}
Let $f\in W(\Omega)$ and $(K_n)_{n\in\mathbb{N}}$ as in
Definition \ref{def-specificacao-determinada-Phi}.
Given $g\in C(\Omega)$ and  $\varepsilon>0$ there is
$n_0\equiv n_0(f,g)\in\mathbb{N}$ such that if
$n\geq n_0$ then
	\[
		\sup_{ y,z\in \Omega }
		\left|\int_{\Omega}g(x)\, dK_{n}(x,y)-\int_{\Omega}g(x)\, dK_{n}(x,z)\right|
		=\mathcal{O}(\varepsilon).
	\]	
 \end{lemma}
\begin{proof}
Given $\varepsilon>0$, follows from the Walters condition \eqref{walters-condition}
that there is $n_1\in\mathbb{N}$ so that if $n\geq n_1$, then
$
	|S_n(f)(x_1^ny_{n+1}^{\infty})-S_n(f)(x_1^nz_{n+1}^{\infty})|
	\leq
	\log(1+\varepsilon),
$
for all $x,y$ and $z\in\Omega$.
Therefore
\[
-\log(1+\varepsilon)
\leq	
S_n(f)(x_1^ny_{n+1}^{\infty})-S_n(f)(x_1^nz_{n+1}^{\infty})
\leq
\log(1+\varepsilon)
\]
which implies that
\[
(1+\varepsilon)^{-1}
\leq	
\frac{\exp(S_n(f)(x_1^ny_{n+1}^{\infty}))}{\exp(S_n(f)(x_1^nz_{n+1}^{\infty}))}
\leq
1+\varepsilon.
\]
From the above inequality follows that
$(1+\varepsilon)^{-1}Z_{n}^{z}\leq Z_{n}^y\leq (1+\varepsilon)Z_{n}^{z}$.
Since $g$ is a continuous function and its domain $\Omega$ is a compact set follows that
$g$ is uniformly continuous, and so there is $n_2\in\mathbb{N}$ such that
if $n\geq n_2$ then
$|g(x_1^nz_{n+1}^{\infty})-g(x_1^ny_{n+1}^{\infty})|<\varepsilon$,
for all $x,y$ and $z\in\Omega$. For all $n\geq n_0\equiv \max\{n_1,n_2\}$
we have
	\begin{align*}
	\int_{\Omega} g(x)\, dK_n(x,z)
	=&
	\frac{1}{Z_{n}^{z}}
		\sum_{ \substack{ x\in\Omega;\\ \sigma^n(x)=\sigma^n(z) }   }
		g(x)\exp(H_n(x))		
	\\
	\leq&
	\frac{(1+\varepsilon)^2}{Z_{n}^{y}}
		\sum_{ \substack{ x\in\Omega;\\ \sigma^n(x)=\sigma^n(y) }   }
		(g(x)+\varepsilon)\exp(H_n(x))		
	\\
	=&
	(1+\varepsilon)^2\int_{\Omega} g(x)\, dK_n(x,y) +(1+\varepsilon)^2\varepsilon
	\\
	=&
	\int_{\Omega} g(x)\, dK_n(x,y) +\mathcal{O}(\varepsilon).
	\end{align*}
By a similar reasoning we obtain the reverse inequality.
\end{proof}

\begin{corollary}\label{cor-aux-unicidade}
Let $f\in W(\Omega)$ and $(K_n)_{n\in\mathbb{N}}$ as in
Definition \ref{def-specificacao-determinada-Phi}.
If $(y_n)_{n\in\mathbb{N}}$ is a sequence in $\Omega$
such that $y_n\to y^*$ and $K_n(\cdot,y_n)\rightharpoonup \tilde{\nu}$,
then $K_n(\cdot,y^{*})\rightharpoonup \tilde{\nu}$.
 \end{corollary}

\begin{proof}
For any fixed $g\in C(\Omega)$ we have
\begin{align*}
\left|\int_{\Omega}g(x)\, dK_{n}(x,y^*)-\int_{\Omega}g(x)\, d\tilde{\nu}(x)\right|
&\leq
\left|\int_{\Omega}g(x)\, dK_{n}(x,y^*)-\int_{\Omega}g(x)\, dK_{n}(x,y_n)\right|
\\
&\qquad +
\left|\int_{\Omega}g(x)\, dK_{n}(x,y_n)-\int_{\Omega}g(x)\, d\tilde{\nu}(x)\right|.
\end{align*}
Given $\varepsilon>0$ follows from Lemma \ref{lemma-aux-unicidade}
that the first term in rhs above is smaller than $\varepsilon$ if
$n$ is large enough. The second term can also be made smaller than
$\varepsilon$ since $K_{n}(\cdot,y_n) \rightharpoonup \tilde{\nu}$.
\end{proof}

\begin{theoremalpha}\label{Teo-DLR-Lim-Term}
Let $f$ be a continuous potential and $(K_n)_{n\in\mathbb{N}}$ as in
Definition \ref{def-specificacao-determinada-Phi}.
Then $\mathcal{G}^{DLR}(f)= \mathcal{G}^{TL}(f)$.
\end{theoremalpha}

\begin{proof}
The inclusion $\mathcal{G}^{TL}(f)\subset\mathcal{G}^{DLR}(f)$
is the content of Proposition \ref{prop-GTL-in-GDLR}.
Suppose by contradiction that there exists
$\mu\in \mathcal{G}^{DLR}(f)$ which is not in $\mathcal{G}^{TL}(f).$
By using the compactness of $\mathcal{G}^{DLR}(f)$ and the classical
hyperplane separation theorem we can ensure the existence of a continuous
function $g:\Omega\to\mathbb{R}$ and $\epsilon>0$ such that
\[
\int_{\Omega}g \, d\mu < \int_{\Omega}g \, d\nu -\epsilon,
\quad \forall \nu\in \mathcal{G}^{TL}(f).
\]
From Theorem \ref{teo-DLR-equation-equiv-gibbs},
for any $n\in\mathbb{N}$, we have
\[
	\int_{\Omega}
		\left[
			\int_{\Omega} g(x) \, dK_{n}(x,y)
		\right]
	\, d\mu(y)
	=
	\int_{\Omega}g\, d\mu.
\]
Therefore, for each $n\in\mathbb{N}$, we have from the previous inequality that
there is $y_n\in \Omega$ such that
\[
\int_{\Omega} g(x) \, dK_{n}(x,y_n)
<
\int_{\Omega} g\, d\nu - \epsilon.
\]
Up to subsequences, we can suppose that
$K_n(\cdot,y_n)\rightharpoonup \tilde{\nu}$
and $y_n\to y^*$. From Corollary \ref{cor-aux-unicidade}
follows that $K_{n}(\cdot,y^*)\rightharpoonup \tilde{\nu}$ and
consequently $\tilde{\nu}\in \mathcal{G}^{TL}(f)$
which is contradiction, thus
showing that $\mathcal{G}^{DLR}(f)= \mathcal{G}^{TL}(f)$.
\end{proof}

\begin{theoremalpha}\label{teo-principal}
If $f\in W(\Omega)$ then
$
\mathcal{G}^{TL}(f)=\mathcal{G}^{DLR}(f)=\mathcal{G}^{*}(f).
$
\end{theoremalpha}

\begin{proof}
If $f\in W(\Omega)$ then we know that $\mathcal{G}^{*}(f)$ is a singleton
(\cite{MR2342978}),
$\# \mathcal{G}^{TL}(f)\geq 1$ and $\mathcal{G}^{TL}(f)=\mathcal{G}^{DLR}(f)$
(Theorem \ref{Teo-DLR-Lim-Term})
so it is enough to prove that
$\mathcal{G}^{TL}(f)\subset\mathcal{G}^{*}(f)$.
From Proposition \ref{prop-Kn-Ln} we have
for any $g\in C(\Omega)$ and $y\in\Omega$ fixed
	\[
		\frac{\mathcal{L}^{n}_{f}(g)(\sigma^n(y))}
		{\mathcal{L}^{n}_{f}(1)(\sigma^n(y))}
		=
		\int_{\Omega} g(x)\, dK_{n}(x,y).
	\]
Assume $K_n$ converges, up to a subsequence, to some probability measure $\nu$.
Then
	\[
	\int_{\Omega} g\, d\nu
	=
	\lim_{n\to\infty}
	\int_{\Omega} g(x)\, dK_{n}(x,y)
	=
	\lim_{n\to\infty}
		\frac{\mathcal{L}^{n}_{f}(g)(\sigma^n(y))}
		{\mathcal{L}^{n}_{f}(1)(\sigma^n(y))}
	=
	\int_{\Omega} g\,  d\nu_{f},
	\]
where the above limit is computed in \cite{MR2342978} and
$\nu_{f}\in \mathcal{G}^{*}(f)$.
Since the function $g\in C(\Omega)$ in above equation is arbitrary,
follows that $\nu=\nu_{f}$, thus finishing the proof.
\end{proof}

\section{Ising Model and Walters Condition} \label{Fro2}

In this section we briefly discuss the long-range Ising model
in Thermodynamic Formalism setting and apply the above results
to ensure the uniqueness of the DLR-Gibbs measures of this model.

The long-range Ising model on the lattice $\mathbb{N}$
with $1/r^{2+\varepsilon}$ interaction energy,
is usually defined by the means of the interaction $\Phi$  of Example
\ref{exemplo-ising-longo-alcance} with $\alpha=2+\varepsilon$, i.e.,
	\[
		\Phi_{A}(x) =
		\begin{cases}
			\displaystyle\frac{x_{n}x_{m}}{|n-m|^{2+\varepsilon}},
					&\ \text{if}\ \ A=\{n,m\}\ \text{and}\ m\neq n;
			\\
			0,			&\ \text{otherwise}.
		\end{cases}
	\]
A straightforward computation shows that the potential
$f:\Omega\to\mathbb{R}$ given by
	\[
		f(x)= \sum_{n\geq 2} \frac{x_1x_{n}}{(n-1)^{2+\varepsilon}}
	\]
is according to Lemma \ref{lema-fcont-interacao} the potential corresponding to $\Phi$.
It is simple to show that $f$ is not $\alpha$-H\"older
continuous for any $0<\alpha<1$. On the other hand, we have
that $f$ is in the Walters class for any $\varepsilon>0$.
Indeed, it is easy to see that for $n,p\in\mathbb{N}$ we
have
\[
\mathrm{var}_{n+p} (S_n(f))
=
(n+p)^{-2-\varepsilon+1} + (n+p-1)^{-2-\varepsilon+1}+...+ p^{-2-\varepsilon+1}.
\]
Therefore, for $p$ fixed, we have
\begin{align*}
\sup_{n\in \mathbb{N}} \
\mathrm{var}_{n+p}(S_n(f))
\leq \mathrm{const.}
\ \sum_{j=p}^\infty j^{-2-\varepsilon+1}
\leq \mathrm{const.} \
p^{-\varepsilon},
\end{align*}
which proves that the potential $f$ is in Walters space.

\medskip

Now, we can apply Theorem \ref{teo-principal}
to ensure that this Ising model has a unique DLR-Gibbs measure
and therefore it has no phase transition in the sense
of multiples DLR-Gibbs measures.

\section{Concluding Remarks}

In this paper we have compared the definitions of Gibbs measures
in terms of the Ruelle operator and specifications.
We show how to obtain for potentials in the Walters and H\"older class
the Gibbs measures usually considered in the Thermodynamic Formalism
via the DLR formalism and prove that the measures obtained from
both approaches are the same.

Both approaches have their advantages.
For example, using the Ruelle operator we were able to prove
some uniform convergence theorems for $K_n(\cdot,y)$.

The literature about absolutely uniformly summable interactions
is vast and this approach allow us to consider non translation invariant
potentials and other lattices than $\mathbb{N}$.
We also show that the long-range Ising model on $\mathbb{N}$
can be studied using the Ruelle operator,
at least when the interaction energy is of the form $1/r^{\alpha}$ with $\alpha>2$.
In these cases, we have proved that the unique Gibbs measure of this model satisfies
$\mathcal{G}^{DLR}(\Phi)=\mathcal{G}^{TL}(\Phi)= \mathcal{G}^*(f)$,
but on the other hand, it is not clear how to treat the cases $1<\alpha\leq 2$
by using the Ruelle operator and what kind of information is
obtainable through this approach.
It is worth pointing out that
treating this model with the DLR approach
is fairly standard, so the connection made here
suggests that more understanding of the DLR Specification theory
can shed light on more general spaces where one
can efficiently use the Ruelle Operator.
Another important feature of the DLR-measure Theory is that it is also
readily appliable to standard Borel spaces, which includes compact and non-compact spaces
\cite{MR2807681}. Some of the results obtained here
can be extended to compact and metric alphabets,
but measurability issues have to be taken into account and
some theorems requires different approach, although the main ideas are contained here,
see \cite{CL-rcontinuas-2016}.

\section*{Acknowledgments}
We would like to express our thanks to
Aernout van Enter, Anthony Quas, Jos\'e Siqueira
and Rodrigo Bissacot for their suggestions and comments on
an early version of this paper. The authors thanks CNPq and FEMAT
for the financial support.



\end{document}